\documentclass[a4paper,12pt]{amsart}

\usepackage[english]{babel}
\usepackage[latin1]{inputenc}
\usepackage{amsmath}
\usepackage{amsthm}
\usepackage{amsfonts}
\usepackage{amssymb, latexsym, graphics, showidx}
\usepackage{color}

\numberwithin{equation}{section}

\newtheorem{thm}{Theorem}[section]

\newtheorem{prop}[thm]{Proposition}
\newtheorem{lem}[thm]{Lemma}

\newcommand{\finedim}{{\unskip\nobreak\hfil\penalty50
   \hskip2em\hbox{}\nobreak\hfil\mbox{$\Box$ \qquad}
   \parfillskip=0pt \finalhyphendemerits=0\par\medskip}}
\newcommand{\R}{\mathbb{R}}
\newcommand{\N}{\mathbb{N}}
\newcommand{\Z}{\mathbb{Z}}

\newcommand{\lam}{\lambda}
\newcommand{\al}{\alpha}

\newcommand{\p}{\partial}

\newcommand{\ep}{\epsilon}
\newcommand{\I}{\mathcal{I}_1}

\newcommand{\beq}{\begin{equation}}
\newcommand{\eeq}{\end{equation}}
\newcommand{\beqs}{\begin{equation*}}
\newcommand{\eeqs}{\end{equation*}}
\newcommand{\beqa}{\begin{eqnarray}}
\newcommand{\eeqa}{\end{eqnarray}}
\newcommand{\beqas}{\begin{eqnarray*}}
\newcommand{\eeqas}{\end{eqnarray*}}

\title[]{Derivation of Orowan's law\\ from the Peierls-Nabarro model}
\email{monneau@cermics.enpc.fr} \email{spatrizi@math.ist.utl.pt}

\begin{document}


\maketitle \centerline{\scshape R\'{e}gis Monneau }
\medskip
{\footnotesize
 \centerline{Universit\'e Paris-Est, CERMICS, Ecole des Ponts ParisTech,}
   \centerline{6-8 avenue Blaise Pascal, Cit\'{e}
Descartes, Champs sur Marne,}
   \centerline{77455 Marne la Vallée Cedex 2,
France}
} 

\medskip

\centerline{\scshape Stefania Patrizi}
\medskip
{\footnotesize
\centerline{Instituto Superior T\'{e}cnico,  Dep. de Matem\'{a}tica}
   \centerline{Av. Rovisco Pais Lisboa, Portugal}

}
\bigskip

\begin{abstract}
In this paper we consider the time dependent Peierls-Nabarro model in dimension one.
This model is a semi-linear integro-differential equation associated to the half Laplacian.
This model describes the evolution of phase transitions associated to dislocations.
At large scale with well separated dislocations, 
we show that the dislocations move at a velocity proportional to the effective stress. 
This implies Orowan's law which claims that the plastic strain velocity 
is proportional to the product of the density of dislocations by the effective stress.
\end{abstract}

\section{Introduction}

\subsection{Setting of the problem}

In this paper we consider a one-dimensional Peierls-Nabarro model,
describing the motion of dislocations in  crystals. In this model dislocations can be seen as phase transitions
of a function $u^\ep$ solving the following equation for $\ep =1$
\begin{equation}\label{uep}
\begin{cases}
\p_{t}
u^\epsilon=\I[u^\epsilon(t,\cdot)]-W'\left(\frac{u^\epsilon}{\epsilon}\right)&\text{in}\quad \R^+\times\R\\
u^\epsilon(0,x)=u_0(x)& \text{on}\quad \R.
\end{cases}
\end{equation}
Here $\I = -(-\Delta)^{\frac12}$ is the half Laplacian whose expression will be made precise later in (\ref{levy})
and $W$ is a one periodic potential which describes the misfit of atoms in the crystal 
created by the presence of dislocations. 
Equation \eqref{uep}  models the dynamics of   
 parallel straight edge dislocation lines in the same slip plane with the same Burgers vector, 
moving with self-interactions. In other words equation (\ref{uep}) simply describes the motion of dislocations 
by relaxation of the total energy (elastic + misfit). 
For a physical introduction to the Peierls-Nabarro
model, see for instance \cite{hl}, \cite{WXM}; we also refer the reader to the paper of Nabarro
\cite{N} which presents an historical tour on the Peierls-Nabarro model.
The Peierls-Nabarro model has been originally introduced as a variational (stationary) model (see \cite{N}).
The model considered in the present paper, i.e. 
the time evolution Peierls-Nabarro model as a gradient flow dynamics has only been introduced 
quite recently, see for instance \cite{MBW} and \cite{Denoual}, and \cite{mp1} where this model is also presented.
See also the paper \cite{dmopt} that initiated several other works about jump-diffusion reaction equations.

In \cite{mp1} we study the limit as $\ep\rightarrow0$ of the viscosity
solution $u^\ep$ of (\ref{uep}) in higher dimensions and with additional periodic terms.
Under certain assumptions, we show in particular that $u^\ep$ converges to the solution of the following equation:
\begin{equation}\label{ueffett}
\begin{cases}
\p_{t} u=\overline{H}(u_x,\I[u(t,\cdot)])
&\text{in}\quad \R^+\times\R\\
u(0,x)=u_0(x)& \text{on}\quad \R.
\end{cases}
\end{equation} 
In mechanics, equation (\ref{ueffett}) can be interpreted as a plastic flow rule, which expresses the plastic strain velocity $\partial_t u$ as a function $\overline{H}$ of the dislocation density $u_x$ and the effective stress $\I[u]$ created by the density of dislocations.
Mathematically the function $\overline{H}$, usually called {\em
effective Hamiltonian}, is determined by the following auxiliary problem:
\begin{equation}\label{sv}
\begin{cases}
\p_{\tau} v= L+\I[v(\tau,\cdot)]-W'(v)&\text{in}\quad \R^+\times\R\\ 
v(0,y)=py& \text{on}\quad \R.
\end{cases}
\end{equation}
Here the quantity $L$ appears to be an additional constant stress field.
Indeed, we have
{ \begin{thm}[Theorem 1.1, \cite{mp1}]\label{ergodic}Assume that $W\in C^{1,1}(\R)$ and $W$ is $1$-periodic. 
For every $L\in\R$ and $p\in\R$, there
exists a unique viscosity solution $v\in C(\R^+\times\R)$ of
\eqref{sv} and there exists a unique $\lam\in\R$ such that $v-py-\lam\tau$
is bounded in $\R^+\times\R$. The real
number  $\lam$ is denoted by $\overline{H}(p,L)$. The function
$\overline{H}(p,L)$ is continuous on $\R^2$ and
non-decreasing in $L$.
\end{thm}}
This is the starting point of this paper.
Our goal is to study the behaviour of $\overline{H}(p,L)$ for small $p$ and $L$,
and in this regime to recover Orowan's law, which claims that
\begin{equation}\label{or}
\overline{H}(p,L) \simeq c_0 |p| L
\end{equation}
for some constant of proportionality $c_0>0$.

\subsection{Main result}

In order to describe our main result, we need the following  assumptions on the potential $W$:
\begin{equation}\label{W}
\begin{cases}W\in C^{4,\beta}(\R)& \text{for some }0<\beta<1\\
W(v+1)=W(v)& \text{for any } v\in\R\\
W=0& \text{on }\Z\\
W>0 & \text{on }\R\setminus\Z\\
\al=W''(0)>0.
\end{cases}
\end{equation}
Under  (\ref{W}), it is in particular known (see Cabr\'{e} and Sol\`{a}-Morales \cite{csm})
that there exists a unique function $\phi$ solution of
\begin{equation}\label{phi}
\begin{cases}\I[\phi]=W'(\phi)&\text{in}\quad \R\\
\phi'>0&\text{in}\quad \R\\
\lim_{x\rightarrow-\infty}\phi(x)=0,\quad\lim_{x\rightarrow+\infty}\phi(x)=1,\quad\phi(0)=\frac{1}{2}.
\end{cases}
\end{equation}
Our main result is the following:
\begin{thm}[Orowan's law]\label{hullprop}Assume \eqref{W} 
and let $p_0,\,L_0\in\R$. Then the function $\overline{H}$ defined in Theorem \ref{ergodic} satisfies
\begin{equation}\label{orowan}\frac{\overline{H}(\delta p_0,\delta L_0)}{\delta^2}\rightarrow
c_0|p_0|L_0\quad\text{as }\delta\rightarrow0^+ \quad \mbox{with}\quad c_0=\left(\int_\R (\phi')^2\right)^{-1}.\end{equation}
\end{thm}
Theorem \ref{hullprop} shows that in the limit of small density of dislocations $p$ and small stress $L$,
the effective Hamiltonian $\overline{H}$ follows Orowan's law (\ref{or}).
This implies that in this regime, the plastic strain velocity $\partial_t u$ in (\ref{ueffett})
is proportional to the dislocation density $|u_x|$ times the effective stress $\I[u]$, i.e.
$$\partial_t u \simeq c_0 |u_x| \I[u(t,\cdot)].$$
Notice that this last equation has been proposed by Head \cite{H} and self-similar solutions have been 
studied mathematically in \cite{bkm}.

Notice that in homogenization problems the effective Hamiltonian is usually
unknown. Explicit formulas for $\overline{H}$ are known only in very special cases, 
see for instance \cite{lpv}. 
The result of Theorem \ref{hullprop} provides an other example of explicit expression 
for a particular homogenization problem.

Finally we give the precise expression (the L\'{e}vy-Khintchine formula in Thm 1 of \cite{di})
of the L\'{e}vy operator $\I$ of order $1$.
For bounded $C^2$- functions $U$ and for $r>0$, we set
\beq\label{levy}\begin{split}\I[U](x) = \I[U,x] &=\int_{|z|\leq
r}(U(x+z)-U(x)-\frac{dU(x)}{dx}\cdot
z)\mu(dz)\\&
+\int_{|z|>r}(U(x+z)-U(x))\mu(dz),\quad \mbox{with}\quad \mu(dz)=\frac1{\pi}\frac{dz}{z^2}.
\end{split}\eeq 
Notice that this expression is independent on the choice of $r>0$, because of the antisymmetry of $z\mu(dz)$.
More generally, when $U$ is $C^2$ such that $U-\ell$ is bounded with $\ell$ a linear function, we simply define
$$ \I[U](x) =\I[U,x]=\lim_{r\to 0^+} \int_{r< |z|< 1/r}(U(x+z)-U(x))\mu(dz)$$

\subsection{Organization of the article}{\mbox{ }}\bigskip

In Section \ref{s-s2}, we present the main ideas which allow us to prove Orowan's law and 
give the proof of the main theorem
(Theorem \ref{hullprop}). This proof is based on Proposition \ref{hproperties} 
which claims asymptotics satisfied by a good Ansatz (see (\ref{eq::s-s100})). 
The remaining part of the paper is then devoted to the proof of Proposition \ref{hproperties}.
In Section \ref{s-s3}, we recall in Lemmata \ref{phiinfinitylem} 
and \ref{psiinfinitylem}, useful asymptotics respectively on the transition layer 
$\phi$ and some corrector $\psi$.
The main result of this section is some asymptotics  on the non linear PDE evaluated on the Ansatz.
In Section \ref{s-s4}, we do the proof of Proposition \ref{hproperties}.
Finally in an appendix (Section \ref{appendix}),
we give the proof of Lemmata \ref{phiinfinitylem} 
and \ref{psiinfinitylem}. We also give the proof of five claims used in Section \ref{s-s3} and a technical lemma (Lemma \ref{lem::s105}) used in Section \ref{s-s4}.

\section{Ideas and  proof of Orowan's law (Theorem \ref{hullprop})}\label{s-s2}{\mbox{ }}

\subsection{Heuristic for the proof of Orowan's law}\label{s2.1}{\mbox{ }}\\

The idea underlying the proof of Orowan's law is related to a fine asymptotics of equation (\ref{sv}). 
It is also known (see \cite{gonzalezmonneau}) that if 
$v$ solves (\ref{sv}) with $L=\delta L_0$, i.e. 
\begin{equation}\label{eq::s102}
\partial_\tau v = \delta L_0 +\I[v(\tau,\cdot)] -W'(v)
\end{equation}
for a choice of initial data with a finite number of indices $i$:
$$\displaystyle v(0,y)=\frac{\delta L_0}{\alpha} + \sum_{x_i^0\ge 0} \phi\left(y-\frac{x_i^0}{\delta}\right)
+ \sum_{x_i^0 < 0} \left(\phi\left(y-\frac{x_i^0}{\delta}\right)-1\right)$$
where $\al=W''(0)>0$ (defined in (\ref{W})), then
$$v^\delta(t,x)=v\left(\frac{t}{\delta^2},\frac{x}{\delta}\right) \to v^0(t,x)
=\sum_{x_i^0\ge 0} H(x-x_i(t)) + \sum_{x_i^0< 0} \left(H(x-x_i(t))-1\right)
\quad \mbox{as}\quad \delta\to 0$$
where $H$ is the Heaviside function and with the dynamics
\begin{equation}\label{eq::s101}
\left\{\begin{array}{l}
\displaystyle \frac{dx_i}{dt}= c_0\left(-L_0 +\frac{1}{\pi}\sum_{j\not= i}\frac{1}{x_i-x_j}\right)\\
\\
x_i(0)=x_i^0.
\end{array}\right.
\end{equation}
Moreover for the choice $p=\delta p_0$ with $p_0 >0$ and $x_i^0=i/p_0$ that we extend formally for all $i\in \Z$, 
we see (at least formally) that 
$$|v(0,y)- \delta p_0 y|\le C_\delta.$$
This suggests also that the infinite sum in (\ref{eq::s101}) should vanish (by antisymmetry) 
and then the mean velocity should be 
$$\frac{dx_i}{dt}\simeq -c_0 L_0$$
i.e., after scaling back
$$v(\tau,y)\simeq \delta p_0 (y-c_1\tau) + \mbox{bounded}$$
with the velocity
$$c_1= \frac{d(x_i/\delta)}{d(t/\delta^2)} \simeq -c_0 L_0 \delta$$
i.e.,
$$v(\tau,y)\simeq \delta p_0 y+\lambda \tau+ \mbox{bounded}\quad \mbox{with}\quad \lambda \simeq \delta^2 c_0 p_0 L_0.$$
We deduce that we should have
$$\frac{v(\tau,y)}{\tau} \to \lambda \simeq \delta^2 c_0 p_0 L_0 \quad \mbox{as}\quad \tau\to +\infty.$$
We see that this $\lambda=\overline{H}(\delta p_0,\delta L_0)$ is exactly the one we expect asymptotically 
in Theorem \ref{hullprop} when $p_0>0$.

\subsection{The ansatz used in the proofs}{\mbox{ }}\\

In the spirit of \cite{fim2}, one may expect to find particular solutions $v$ of (\ref{eq::s102}) that we can write
$$v(\tau,y)=h(\delta p_0 y + \lambda \tau)$$
for some $\lambda\in\R$ and 
a function $h$ (called hull function) satisfying
$$|h(z)-z|\le C.$$
This means that $h$ solves
$$\lambda h' = \delta L_0+ \delta |p_0| \I [h] -W'(h).$$
Then it is natural to introduce the non linear operator:
\begin{equation}\label{eq::s103}
NL^\lambda_{L_0}[h]:= \lambda h' - \delta L_0- \delta |p_0| \I [h] + W'(h)
\end{equation}
and for the ansatz for $\lambda$:
$$\overline{\lambda}^{L_0}_\delta = \delta^2 c_0 |p_0|L_0$$
it is natural to look for an ansatz $h^{L_0}_\delta$ for $h$.
The answer is indicated by the heuristic of subsection \ref{s2.1}.
Indeed we define (see Proposition \ref{hproperties})
$$h_\delta^{L_0}(x)=\lim_{n\rightarrow+\infty}s_{\delta,n}^{L_0}(x)$$
where for all
$p_0\neq 0$, $L_0\in\R$, $\delta>0$ and $n\in\N$ we
define the sequence of functions $\{s_{\delta,n}^{L_0}(x)\}_n$ by
\begin{equation}\label{eq::s-s100}
s_{\delta,n}^{L_0}(x)=\frac{\delta L_0}{\al}+\sum_{i=-n}^{n}\left[\phi\left(\frac{x-i}{\delta|p_0|}\right)
+\delta\psi\left(\frac{x-i}{\delta|p_0|}\right)\right]-n
\end{equation}
where $\al=W''(0)>0$, $\phi$ is the solution of \eqref{phi} and 
the corrector $\psi$ is the solution of the following
problem
\begin{equation}\label{psi}
\begin{cases}\I[\psi]=W''(\phi)\psi+\frac{L_0}{W''(0)}(W''(\phi)-W''(0))+c\phi'&\text{in}\quad \R\\
\lim_{x\rightarrow{+\atop -}\infty}\psi(x)=0\\
c=\frac{L_0}{\int_{\R} (\phi')^2}.
\end{cases}
\end{equation} 
From \cite{gonzalezmonneau}, it is known that there exists a unique $\psi$ solution of (\ref{psi}).
Moreover this corrector $\psi$ has been introduced naturally in \cite{gonzalezmonneau} 
in order to perform part of the analysis presented in the heuristic (subsection \ref{s2.1}),
and this is then natural to use it here in our ansatz.
We will prove later the following result which justifies that the ansatz is indeed a good ansatz as expected.
 \begin{prop}{ \bf (Good ansatz)}\label{hproperties}\\
Assume \eqref{W}. For any $L\in\R$, $\delta>0$ and $x\in\R$, there exists
the finite limit
$$h_\delta^{L}(x)=\lim_{n\rightarrow+\infty}s_{\delta,n}^{L}(x).$$
Moreover $h_\delta^{L}$ has the following properties:
\begin{itemize}
    \item [(i)]$h_\delta^{L}\in C^{2}(\R)$ and
satisfies
$$NL_{L}^{\overline{\lam}_\delta^{L}}[h_\delta^{L}](x)=o(\delta),$$ 
where
 $\lim_{\delta\rightarrow 0}\frac{o(\delta)}{\delta}=0$, uniformly
 for $x\in\R$ and locally uniformly in $L\in\R$;
Here 
$$\overline{\lam}_\delta^{L} = \delta^2 c_0 |p_0| L $$
and $NL_{L}^{\lambda}$ is defined in (\ref{eq::s103}).
    \item [(ii)]There exists a constant $C>0$ such that
    $|h_\delta^{L}(x)-x|\leq C$ for any $x\in\R$.
\end{itemize}
\end{prop}




\subsection{Proof of Theorem \ref{hullprop}}

We will show that Theorem \ref{hullprop} follows from Proposition \ref{hproperties},
and the comparison principle.

Fix $\eta>0$ and let $L=L_0-\eta$. By (i) of Proposition
\ref{hproperties}, there exists $\delta_0=\delta_0(\eta)>0$ such
that for any $\delta\in(0, \delta_0)$ we have
\begin{equation}\label{NLhL_0}
NL_{L_0}^{\overline{\lam}_\delta^{L}}[h_\delta^L]=NL_{L}^{\overline{\lam}_\delta^{L}}[h_\delta^L]-\delta\eta<0\quad\text{in
}\R.
\end{equation}
Let us consider the function $\widetilde{v}(\tau,y)$, defined by
$$\widetilde{v}(\tau,y) = h_\delta^L(\delta p_0 y + \overline{\lambda}^L_\delta \tau).$$
By (ii) of Proposition \ref{hproperties}, we have
\begin{equation}\label{propubarra}|\widetilde{v}(\tau,y)-\delta p_0 y - \overline{\lam}_\delta^L\tau|\leq \lceil
C\rceil,\end{equation} 
{ where $\lceil C\rceil$ is the ceil integer part of $C$.}
Moreover, by \eqref{NLhL_0} and \eqref{propubarra}, $\widetilde{v}$
satisfies
$$\left\{%
\begin{array}{ll}
   \widetilde{v}_\tau\leq \delta L_0 + \I[\widetilde{v}]-W'(\widetilde{v})& \hbox{in } \R^+\times\R\\
    \widetilde{v}(0,y) \leq \delta p_0 y+ \lceil C\rceil & \hbox{on } \R.\\
\end{array}%
\right.$$
 Let $v(\tau,y)$ be the solution of \eqref{sv}, { with $p=\delta p_0$ and $L=\delta L_0$}, whose existence is ensured by Theorem
\ref{ergodic}.
 Then from the comparison principle and the
periodicity of $W$, we deduce that
$$\widetilde{v}(\tau,y)\leq v(\tau,y)+\lceil C\rceil.$$
By the previous inequality and \eqref{propubarra}, we
get
$$\overline{\lam}_\delta^L\tau\leq v(\tau,y)-\delta p_0 y +2\lceil
C\rceil,$$ and dividing by $\tau$ and letting $\tau$ go to
$+\infty$, we finally obtain
$$\delta^2c_0|p_0|(L_0-\eta)=\overline{\lam}_\delta^L\leq \overline{H}(\delta p_0,\delta
L_0).$$Similarly, it is possible to show that
$$\overline{H}(\delta p_0,\delta
L_0)\leq \delta^2c_0|p_0|(L_0+\eta).$$We have proved that for any
$\eta>0$ there exists $\delta_0=\delta_0(\eta)>0$ such that for
any $\delta\in(0,\delta_0)$ we have
$$\left|\frac{\overline{H}(\delta p_0,\delta
L_0)}{\delta^2}-c_0|p_0|L_0\right|\leq c_0|p_0|\eta,$$ i.e.
\eqref{orowan}, as desired.\finedim

\section{Preliminary asymptotics}\label{s-s3}

The main goal of this section is to show Lemma \ref{lem::1011}
which is a first result in the direction of Proposition \ref{hproperties}.
We start with prelimary results in a first subsection and prove Lemma \ref{lem::1011}
in the second subsection.

\subsection{Preliminary results}\mbox{ } \bigskip

On the function $W$, we assume
(\ref{W}). Then there exists a unique solution of \eqref{phi}
which is of class  $C^{2,\beta}$, as shown by Cabr\'{e} and
Sol\`{a}-Morales in \cite{csm}. Under \eqref{W}, the existence of
a solution of class $C^{1,\beta}$ of the problem \eqref{psi} is
proved by Gonz\'{a}les and Monneau in \cite{gonzalezmonneau}.
Actually, the regularity of $W$ implies, that $\phi\in
C^{4,\beta}(\R)$ and $\psi\in C^{3,\beta}(\R)$, see Lemma 2.3 in
\cite{csm}.

To prove Proposition \ref{hproperties} we need several preliminary
results. We first state the following two lemmata about the
behavior of the functions $\phi$ and $\psi$ at infinity. We denote
by $H(x)$ the Heaviside function defined by
\begin{equation*}
H(x)=\begin{cases}1 &\text{for }x\geq 0 \\ 0& \text{for
}x<0.\end{cases}
\end{equation*}Then we have
\begin{lem}[Behavior of $\phi$]\label{phiinfinitylem}Assume \eqref{W}. Let $\phi$ be the solution of
\eqref{phi}, then there exist constants $K_0,K_1 >0$ such that
\begin{equation}\label{phiinfinity}\left|\phi(x)-H(x)+\frac{1}{\al \pi
x}\right|\leq \frac{K_1}{x^2},\quad\text{for }|x|\geq 1,
\end{equation}and for any $x\in\R$
\begin{equation}\label{phi'infinity}0<\frac{K_0}{1+x^2}\leq
\phi'(x)\leq\frac{K_1}{1+x^2},\end{equation}
\begin{equation}\label{phi''infinity}-\frac{K_1}{1+x^2}\leq
\phi''(x)\leq\frac{K_1}{1+x^2},\end{equation}
\begin{equation}\label{phi'''infinity}-\frac{K_1}{1+x^2}\leq
\phi'''(x)\leq\frac{K_1}{1+x^2}.\end{equation}
\end{lem}
\begin{lem}[Behavior of $\psi$]\label{psiinfinitylem}Assume \eqref{W}. Let $\psi$ be the solution of
\eqref{psi}, then for any $L\in\R$  there exist constants $K_2$
and $K_3$, with $K_3>0$, depending on $L$ such that
\begin{equation}\label{psiinfinity}\left|\psi(x)-\frac{K_2}{
x}\right|\leq\frac{K_3}{x^2},\quad\text{for }|x|\geq 1,
\end{equation}and for any $x\in\R$
\begin{equation}\label{psi'infinity}-\frac{K_3}{1+x^2}\leq
\psi'(x)\leq \frac{K_3}{1+x^2},
\end{equation}
\begin{equation}\label{psi''infinity}-\frac{K_3}{1+x^2}\leq
\psi''(x)\leq \frac{K_3}{1+x^2}.
\end{equation}
\end{lem}We postpone the proof of the two lemmata in the appendix (Section \ref{appendix}).

For simplicity of notation we denote (for the rest of the paper)
$$x_i=\frac{x-i}{\delta|p_0|},\quad\widetilde{\phi}(z)=\phi(z)-H(z),\quad \I[\phi,x_i]=\I[\phi](x_i).$$
 Then we have the following five claims
(whose proofs are also postponed in the appendix (Section \ref{appendix})).\\
\noindent {\bf Claim 1:} {\em Let $x=i_0+\gamma$, with $i_0\in\Z$
and $\gamma\in\left(-\frac{1}{2},\frac{1}{2}\right]$, then
\begin{equation*}\sum_{i=-n\atop i\neq
i_0}^n\frac{1}{x-i}\rightarrow -2\gamma\sum_{i=1}^{
+\infty}\frac{1}{i^2-\gamma^2}\quad\text{as
}n\rightarrow+\infty,\end{equation*}
\begin{equation*}
\sum_{i=-n\atop }^{i_0-1}\frac{1}{(x-i)^2}\rightarrow
\sum_{i=1}^{+\infty}\frac{1}{(i+\gamma)^2}\quad\text{as
}n\rightarrow+\infty,\end{equation*}
\begin{equation*}\sum_{i=i_0+1\atop }^{n}\frac{1}{(x-i)^2}\rightarrow
\sum_{i=1}^{+\infty}\frac{1}{(i-\gamma)^2}\quad\text{as
}n\rightarrow+\infty.\end{equation*}}\\

\noindent {\bf Claim 2:} {\em For any $x\in\R$ the sequence
$\{s_{\delta,n}^L(x)\}_n$ converges as $n\rightarrow+\infty$.}\\

\noindent {\bf Claim 3:} {\em The sequence
$\{(s_{\delta,n}^L)'\}_n$ converges on $\R$ as
$n\rightarrow+\infty$, uniformly on compact sets.} \\

\noindent {\bf Claim 4:} {\em The sequence
$\{(s_{\delta,n}^L)''\}_n$ converges
on $\R$ as $n\rightarrow+\infty$, uniformly on compact sets.}\\

\noindent {\bf Claim 5:} {\em For any $x\in\R$ the sequences
$\sum_{i=-n}^n\I[\phi,x_i]$ and $\sum_{i=-n}^n\I[\psi,x_i]$
converge as $n\rightarrow+\infty$.}\\

\subsection{First asymptotics}\mbox{ } \bigskip

In order to do the proof of Proposition \ref{hproperties}, we
first get the following result:
\begin{lem}{\bf (First asymptotics)}\label{lem::1011}
We have
$$-C\delta^2\leq\lim_{n\rightarrow+\infty}NL_L^{\overline{\lam}_\delta^L}[s_{\delta,n}^L](x)\leq
C\delta^2,$$ 
where $C$ is independent of $x$.
\end{lem}

\noindent {\bf Proof of Lemma \ref{lem::1011}.}\\

\noindent {\bf Step 1: First computation}\\
Fix $x\in\R$, let $i_0\in\Z$ and
$\gamma\in\left(-\frac{1}{2},\frac{1}{2}\right]$ be such that
$x=i_0+\gamma$, let $\frac{1}{\delta|p_0|}\geq2$ and $n>|i_0|$.
Then we have 

\begin{equation*}\begin{split}
A:= \quad &NL_L^{\overline{\lam}_\delta^L}[s_{\delta,n}^L](x)\\
\\
=\quad &\frac{\overline{\lam}_\delta^L}{\delta|p_0|}\sum_{i=-n}^n\left[\phi'(x_i)+\delta
\psi'(x_i)\right]-\sum_{i=-n}^n\left[\I[\phi,x_i]+\delta\I[\psi,x_i]\right]\\&+W'\left(\frac{L\delta}{\al}
+\sum_{i=-n}^n\left[\phi(x_i)+\delta
\psi(x_i)\right]\right)-\delta L
\end{split}\end{equation*}
where we have used the definitions and the periodicity of $W$.
Using the equation (\ref{phi}) satisfied by $\phi$, we can rewrite it as

\begin{equation*}\begin{split}
A=\quad &\frac{\overline{\lam}_\delta^L}{\delta|p_0|}\left\{\phi'(x_{i_0})+\delta
\psi'(x_{i_0})+\sum_{i=-n\atop i\neq i_0}^n\left[\phi'(x_i)+\delta
\psi'(x_i)\right]\right\} -\sum_{i=-n\atop i\neq
i_0}^nW'(\widetilde{\phi}(x_i))-\delta\I[\psi,x_{i_0}]\\
&-\delta\sum_{i=-n\atop
i\neq i_0}^n\I[\psi,x_{i}]+W'\left(\frac{L\delta}{\al}
+\sum_{i=-n}^n\left[\widetilde{\phi}(x_i)+\delta
\psi(x_i)\right]\right)-W'(\widetilde{\phi}(x_{i_0}))-\delta
L
\end{split}\end{equation*}
Using the definition of $\overline{\lam}_\delta^L$ and a Taylor expansion of $W'$, we get

\begin{equation*}\begin{split}
A=\quad &\delta c_0L\left\{\phi'(x_{i_0})+\delta
\psi'(x_{i_0})+\sum_{i=-n\atop i\neq i_0}^n\left[\phi'(x_i)+\delta
\psi'(x_i)\right]\right\}-W''(0)\sum_{i=-n\atop i\neq
i_0}^n\widetilde{\phi}(x_i)-\delta\sum_{i=-n\atop i\neq
i_0}^n\I[\psi,x_{i}]\\
&-\delta\I[\psi,x_{i_0}]+W''(\phi(x_{i_0}))\left(\frac{L\delta}{\al}
+\delta \psi(x_{i_0})+\sum_{i=-n\atop i\neq
i_0}^n\left[\widetilde{\phi}(x_i)+\delta \psi(x_i)\right]\right)-\delta L +E
\end{split}\end{equation*}

with the error term
$$E= \sum_{i=-n\atop i\neq
i_0}^nO(\widetilde{\phi}(x_i))^2+O\left(\frac{L\delta}{\al}
+\delta \psi(x_{i_0})+\sum_{i=-n\atop i\neq
i_0}^n\left[\widetilde{\phi}(x_i)+\delta
\psi(x_i)\right]\right)^2$$
Simply reorganizing the terms, we get with $c=c_0 L$:

\begin{equation*}\begin{split}
A=\quad &\delta c_0L\left\{\delta
\psi'(x_{i_0})+\sum_{i=-n\atop i\neq i_0}^n\left[\phi'(x_i)+\delta
\psi'(x_i)\right]\right\}-W''(0)\sum_{i=-n\atop i\neq
i_0}^n\widetilde{\phi}(x_i)-\delta\sum_{i=-n\atop i\neq
i_0}^n\I[\psi,x_{i}]\\
&+W''(\phi(x_{i_0}))\left(\sum_{i=-n\atop
i\neq i_0}^n\left[\widetilde{\phi}(x_i)+\delta
\psi(x_i)\right]\right)\\
&+\delta\Big(-\I[\psi,x_{i_0}]+W''(\phi(x_{i_0}))\psi(x_{i_0})+\frac{L}{\al}W''(\phi(x_{i_0}))-L+c\phi'(x_{i_0})\Big)+E
\end{split}\end{equation*}

Using equation (\ref{psi}) satisfied by $\psi$, we get 

\begin{equation*}\begin{split}
A=\quad &\delta c_0L\left\{\delta \psi'(x_{i_0})+\sum_{i=-n\atop i\neq
i_0}^n\left[\phi'(x_i)+\delta
\psi'(x_i)\right]\right\}+(W''(\phi(x_{i_0}))-W''(0))\sum_{i=-n\atop
i\neq i_0}^n\widetilde{\phi}(x_i)\\
&-\delta\sum_{i=-n\atop i\neq
i_0}^n\I[\psi,x_{i}]+W''(\phi(x_{i_0}))\delta\sum_{i=-n\atop i\neq
i_0}^n \psi(x_i)+E
\end{split}\end{equation*}

\noindent {\bf Step 2: Bound on $\sum_{i=-n\atop i\neq i_0}^n\left[\phi'(x_i)+\delta
\psi'(x_i)\right]$}\\

Let us bound the second term  of the last equality, uniformly in
$x$. From \eqref{phi'infinity} and \eqref{psi'infinity} it follows
that
$$-\delta^3|p_0|^2 K_3\sum_{i=-n\atop
i\neq i_0}^n\frac{1}{(x-i)^2}\leq\sum_{i=-n\atop i\neq
i_0}^n\left[\phi'(x_i)+\delta \psi'(x_i)\right]\leq
\delta^2|p_0|^2(K_1+\delta K_3)\sum_{i=-n\atop i\neq
i_0}^n\frac{1}{(x-i)^2},$$ and then by Claim 1 we get
\begin{equation}\label{NL1}-C\delta^3\leq\lim_{n\rightarrow+\infty}\sum_{i=-n\atop i\neq i_0}^n\left[\phi'(x_i)+\delta
\psi'(x_i)\right]\leq C\delta^2.\end{equation} Here and
henceforth, $C$ denotes various positive constants independent of
$x$.\\

\noindent {\bf Step 3: Bound on $(W''(\phi(x_{i_0}))-W''(0))\sum_{i=-n\atop
i\neq i_0}^n\widetilde{\phi}(x_i)$}\\

Now, let us prove that
\begin{equation}\label{NL2}-C\delta^2\leq\lim_{n\rightarrow+\infty}(W''(\phi(x_{i_0}))-W''(0))\sum_{i=-n\atop
i\neq i_0}^n\widetilde{\phi}(x_i)\leq C\delta^2.
\end{equation}
By \eqref{phiinfinity} we have
\begin{equation}\label{nlphi1}\begin{split}\left|\sum_{i=-n\atop i\neq
i_0}^n\widetilde{\phi}(x_i)+\frac{\delta|p_0|}{\al\pi}\sum_{i=-n\atop
i\neq i_0}^n\frac{1}{x-i}\right|\leq
K_1\delta^2|p_0|^2\sum_{i=-n\atop i\neq i_0}^n\frac{1}{(x-i)^2}.
\end{split}\end{equation} If $|\gamma|\geq \delta|p_0|$,  then
again from \eqref{phiinfinity}, $|\widetilde{\phi}(x_{i_0})
+\frac{\delta|p_0|}{\al\pi\gamma}|\leq
K_1\frac{\delta^2|p_0|^2}{\gamma^2}$ which implies that
$$|W''(\widetilde{\phi}(x_{i_0}))-W''(0)|\leq
|W'''(0)\widetilde{\phi}(x_{i_0})|+O(\widetilde{\phi}(x_{i_0}))^2\leq
C\frac{\delta}{|\gamma|}+C\frac{\delta^2}{\gamma^2}.$$ By the
previous inequality, \eqref{nlphi1} and Claim 1 we deduce that
\begin{equation*}\left|\lim_{n\rightarrow+\infty}(W''(\phi(x_{i_0}))-W''(0))\sum_{i=-n\atop
i\neq i_0}^n\widetilde{\phi}(x_i)\right|\leq
C\left(\frac{\delta}{|\gamma|}+\frac{\delta^2}{\gamma^2}\right)(\delta|\gamma|+\delta^2)\leq
C\delta^2,
\end{equation*}where $C$ is independent of $\gamma$.

Finally, if $|\gamma|<\delta|p_0|$, from \eqref{nlphi1} and Claim
1 we conclude that
\begin{equation*}\left|\lim_{n\rightarrow+\infty}(W''(\phi(x_{i_0}))-W''(0))\sum_{i=-n\atop
i\neq i_0}^n\widetilde{\phi}(x_i)\right|\leq
C\delta|\gamma|+C\delta^2\leq C\delta^2,
\end{equation*}and \eqref{NL2} is proved.\\

\noindent {\bf Step 4: Bound on $\delta\sum_{i=-n\atop i\neq
i_0}^n\I[\psi,x_{i}]$}\\
We have
\begin{equation}\label{i1psi}\I[\psi]=W''(\widetilde{\phi})\psi+\frac{L}{\al}(W''(\widetilde{\phi})-W''(0))+c\phi'
=W''(0)\psi+\frac{L}{\al}W'''(0)\widetilde{\phi}+O(\widetilde{\phi})\psi+O(\widetilde{\phi})^2+c\phi'.\end{equation}


Then by \eqref{i1psi}, \eqref{phiinfinity}, \eqref{phi'infinity},
\eqref{psiinfinity} and Claim 1, we have
\begin{equation}\label{NL3}\left|\lim_{n\rightarrow+\infty}\delta\sum_{i=-n\atop i\neq
i_0}^n\I[\psi,x_{i}]\right|\leq C\delta^2.
\end{equation}

\noindent {\bf Step 5: Bound on $W''(\phi(x_{i_0}))\delta\sum_{i=-n\atop
i\neq i_0}^n \psi(x_i)$}\\
Similarly
\begin{equation}\label{NL4}\left|\lim_{n\rightarrow+\infty}W''(\phi(x_{i_0}))\delta\sum_{i=-n\atop
i\neq i_0}^n \psi(x_i)\right|\leq C\delta^2.
\end{equation}
\noindent {\bf Step 6: Bound on the remaining part $E$}\\
Finally, still from \eqref{phiinfinity}, \eqref{psiinfinity},
 and Claim 1 it follows
that
\begin{equation}\label{NL5}\left|\lim_{n\rightarrow+\infty}\sum_{i=-n\atop i\neq
i_0}^n(\widetilde{\phi}(x_i))^2+O\left(\frac{L\delta}{\al} +\delta
\psi(x_{i_0})+\sum_{i=-n\atop i\neq
i_0}^n\left[\widetilde{\phi}(x_i)+\delta
\psi(x_i)\right]\right)^2\right|\leq C\delta^2.
\end{equation}
\noindent {\bf Step 7: Conclusion}\\
Therefore, from \eqref{NL1}, \eqref{NL2}, \eqref{NL3}, \eqref{NL4} and
\eqref{NL5} we conclude that
$$-C\delta^2\leq \lim_{n\rightarrow+\infty}NL_L^{\overline{\lam}_\delta^L}[s_{\delta,n}^L]\leq C\delta^2$$ with $C$ independent of $x$ and
Lemma \ref{lem::1011} is
proved.\\

\section{Proof of Proposition \ref{hproperties}}\label{s-s4}
{\mbox{ }} \bigskip

In order to perform the proof of Proposition \ref{hproperties},
we will use the following technical result whose proof is postponed to the appendix.
\begin{lem}{\bf (Vanishing far away contribution)}\label{lem::s105}\\
We have
\begin{equation}\label{limalimnI2}\lim_{a\rightarrow+\infty}\lim_{n\rightarrow+\infty}\int_{|y|\geq
a}[s_{\delta,n}^L(x+y)-s_{\delta,n}^L(x)]\mu(dy)=0.\end{equation}
\end{lem}

{ 
We also need to introduce the notation
$$\I^1[f,x]= \int_{|y|<1}[f(x+y)-f(x)-f'(x)y]\mu(dy)$$
and
$$\I^2[f,x]= \int_{|y|\geq1}[f(x+y)-f(x)]\mu(dy).$$
}

\noindent {\bf Proof of Proposition \ref{hproperties}}\\
\noindent {\bf Step 1: proof of ii)}\\
Let   $x=i_0+\gamma$ with $i_0\in\Z$ and
$\gamma\in\left(-\frac{1}{2},\frac{1}{2}\right]$. Let $\frac{1}{\delta|p_0|}\geq2$ and $n>|i_0|$,
then by \eqref{phiinfinity} and \eqref{psiinfinity} we get
\begin{equation*}\begin{split}s_{\delta,n}^L(x)-x&=\frac{L\delta}{\al}+\phi(x_{i_0})+\delta\psi(x_{i_0})-n-i_0-\gamma+\sum_{i=-n\atop i\neq
i_0}^{n}[\phi(x_i)
+\delta\psi(x_i)]\\&=\frac{L\delta}{\al}+\phi(x_{i_0})+\delta\psi(x_{i_0})-\gamma+\sum_{i=-n\atop
i\neq i_0}^{n}[\widetilde{\phi}(x_i) +\delta\psi(x_i)]\\&\leq
\frac{L\delta}{\al}+\frac{3}{2}+\delta\|\psi\|_\infty+\sum_{i=-n\atop
i\neq i_0}^{n}\left[-\left(\frac{1}{\al\pi}-\delta
K_2\right)\frac{\delta|p_0|}{x-i}+(K_1+\delta
K_3)\frac{\delta^2|p_0|^2}{(x-i)^2}\right].
\end{split}\end{equation*} Then, by Claim 1
$$h_\delta^L(x)-x=\lim_{n\rightarrow+\infty}s_{\delta,n}^L(x)-x\leq C.$$
Similarly we can prove that $$h_\delta^L(x)-x\geq -C,$$which
concludes the proof of ii).\\

\noindent {\bf Step 2: proof of i)}\\
The function
$h_\delta^L(x)=\lim_{n\rightarrow+\infty}s_{\delta,n}^L(x)$ is
well defined for any $x\in\R$ by Claim 2. Moreover, by Claim 3 and
4 and classical analysis results, it is of class $C^2$ on $\R$
with
$$(h_\delta^{L})'(x)=\lim_{n\rightarrow+\infty}(s_{\delta,n}^L)'(x)=\lim_{n\rightarrow+\infty}
\frac{1}{\delta|p_0|}\sum_{i=-n}^{n}\left[\phi'\left(\frac{x-i}{\delta|p_0|}\right)
+\delta\psi'\left(\frac{x-i}{\delta|p_0|}\right)\right],$$
$$(h_\delta^{L})''(x)=\lim_{n\rightarrow+\infty}(s_{\delta,n}^L)''(x)=\lim_{n\rightarrow+\infty}
\frac{1}{\delta^2|p_0|^2}\sum_{i=-n}^{n}\left[\phi''\left(\frac{x-i}{\delta|p_0|}\right)
+\delta\psi''\left(\frac{x-i}{\delta|p_0|}\right)\right],$$and the
convergence of $\{s_{\delta,n}^L\}_n$, $\{(s_{\delta,n}^L)'\}_n$
and $\{(s_{\delta,n}^L)''\}_n$ is uniform on compact sets.

Let us show that for any $x\in\R$
\begin{equation}\label{NLh}\I[h_\delta^L,x]=\lim_{n\rightarrow+\infty}\I[s_{\delta,n}^L,x].\end{equation}

\noindent {\bf Step 2.1: term $\I^1[h_\delta^L,x]$}\\
First, we prove that
\begin{equation}\label{NLh^1}\I^1[h_\delta^L,x]=\lim_{n\rightarrow+\infty}\I^1[s_{\delta,n}^L,x].\end{equation}
Fix $x\in\R$, we know that for any $y\in[-1,1]$, $y\neq 0$
$$\frac{s_{\delta,n}^L(x+y)-s_{\delta,n}^L(x)-(s_{\delta,n}^L)'(x)y}{|y|^2}\rightarrow
\frac{h_\delta^L(x+y)-h_\delta^L(x)-(h_\delta^{L})'(x)y}{|y|^2}\quad\text{as
}n\rightarrow+\infty.$$  By the uniform convergence of the
sequence
 $\{(s_{\delta,n}^L)''\}_n$ we have
$$\frac{|s_{\delta,n}^L(x+y)-s_{\delta,n}^L(x)-(s_{\delta,n}^L)'(x)y|}{|y|^2}\leq
\sup_{z\in[x-1,x+1]}(s_{\delta,n}^L)''(z)\leq C,$$ where $C$ is
indipendent of $n$, and \eqref{NLh^1} follows from the dominate
convergence Theorem.\\

 \noindent {\bf Step 2.2: term $\I^2[h_\delta^L,x]$}\\
Then, to prove \eqref{NLh} it suffices to show that
\begin{equation*}\I^2[h_\delta^L,x]=\lim_{n\rightarrow+\infty}\I^2[s_{\delta,n}^L,x].\end{equation*}
From Claim 5 and \eqref{NLh^1}, we know that for any $x\in\R$
there exists $\lim_{n\rightarrow+\infty}\I^2[s_{\delta,n}^L,x]$.
For $a>1$, we have
\begin{equation*}\I^2[s_{\delta,n}^L,x]=\int_{1\leq|y|\leq a}[s_{\delta,n}^L(x+y)-s_{\delta,n}^L(x)]\mu(dy)
+\int_{|y|\geq
a}[s_{\delta,n}^L(x+y)-s_{\delta,n}^L(x)]\mu(dy).\end{equation*}
By the uniform convergence of $\{s_{\delta,n}^L\}_n$ on compact
sets
$$\lim_{n\rightarrow+\infty}\int_{1\leq|y|\leq
a}[s_{\delta,n}^L(x+y)-s_{\delta,n}^L(x)]\mu(dy)=\int_{1\leq|y|\leq
a}[h_\delta^L(x+y)-h_\delta^L(x)]\mu(dy),$$then there exists the
limit
$$\lim_{n\rightarrow+\infty}\int_{|y|\geq
a}[s_{\delta,n}^L(x+y)-s_{\delta,n}^L(x)]\mu(dy).$$

 Then, we finally  get
\begin{equation*}\begin{split}\lim_{n\rightarrow+\infty}\I^2[s_{\delta,n}^L,x]
&=\lim_{a\rightarrow+\infty}\lim_{n\rightarrow+\infty}\I^2[s_{\delta,n}^L,x]
\\&=\lim_{a\rightarrow+\infty}\lim_{n\rightarrow+\infty}\int_{1\leq|y|\leq
a}[s_{\delta,n}^L(x+y)-s_{\delta,n}^L(x)]\mu(dy)\\&+\lim_{a\rightarrow+\infty}\lim_{n\rightarrow+\infty}\int_{|y|>
a}[s_{\delta,n}^L(x+y)-s_{\delta,n}^L(x)]\mu(dy)\\&
=\lim_{a\rightarrow+\infty}\int_{1\leq|y|\leq
a}[h_\delta^L(x+y)-h_\delta^L(x)]\mu(dy)\\&=\I^2[h_\delta^L,x],
\end{split}\end{equation*}as
desired, where we have used Lemma \ref{lem::s105}.\\

\noindent {\bf Step 2.3: conclusion}\\
Now we can conclude the proof of (i). Indeed, by Claim 2, Claim 3
and \eqref{NLh}, for any $x\in\R$
$$NL_L^{\overline{\lam}_\delta^L}[h_\delta^L](x)=\lim_{n\rightarrow+\infty}NL_L^{\overline{\lam}_\delta^L}[s_{\delta,n}^L](x),$$
and Lemma \ref{lem::1011} implies that
$$NL_L^{\overline{\lam}_\delta^L}[h_\delta^L](x)=o(\delta),\quad\text{as
}\delta\rightarrow0,$$where
$\lim_{\delta\rightarrow0}\frac{o(\delta)}{\delta}=0$, uniformly
for $x\in\R$.

\section{Appendix}\label{appendix}

In this appendix, we prove the following technical results used in the previous section: Lemmata \ref{phiinfinitylem}
and \ref{psiinfinitylem}, the Claims 1-5 and Lemma \ref{lem::s105}.

\subsection{Proof of Lemma \ref{phiinfinitylem}.}
Properties \eqref{phiinfinity} and \eqref{phi'infinity} are proved
in \cite{gonzalezmonneau}.

Let us show \eqref{phi''infinity}.

 For $a>0$, we denote by
$\phi_a'(x)=\phi'\left(\frac{x}{a}\right)$. Remark that $\phi_a'$
is a solution of
$$\I[\phi_a']=\frac{1}{a}W''(\phi_a)\phi_a'\quad\text{in }\R.$$

Since $\phi''$ is bounded and of class $C^{2,\beta}$, $\I[\phi'']$
is well defined and by deriving twice the equation in \eqref{phi}
we see that $\phi''$ is a solution of
$$\I[\phi'']=W''(\phi)\phi''+W'''(\phi)(\phi')^2.$$ Let
$\overline{\phi}=\phi''-C\phi_a'$, with $C>0$, then
$\overline{\phi}$ satisfies
\beqs\begin{split}\I[\overline{\phi}]-W''(\phi)\overline{\phi}&=C\phi'_a\left(W''(\phi)-\frac{1}{a}W''(\phi_a)\right)+W'''(\phi)(\phi')^2
\\&=C\phi'_a\left(W''(\phi)-\frac{1}{a}W''(\phi_a)\right)+o\left(\frac{1}{1+x^2}\right),\end{split}\eeqs
as $|x|\rightarrow+\infty$, by \eqref{phi'infinity}. Fix $a>0$ and
$R>0$ such that
\begin{equation}\label{ultimalemmaphi}\left\{%
\begin{array}{ll}
    W''(\phi)-\frac{1}{a}W''(\phi_a)>\frac{1}{2}W''(0)>0 & \hbox{on } \R\setminus[-R,R];\\
    W''(\phi)>0, & \hbox{on } \R\setminus[-R,R]. \\
\end{array}%
\right.\end{equation}  Then from \eqref{phi'infinity}, for $C$
large enough we get
$$\I[\overline{\phi}]-W''(\phi)\overline{\phi}\geq 0\quad\text{on
}\R\setminus[-R,R].$$ Choosing $C$ such that moreover
$$\overline{\phi}<0\quad \text{on
}[-R,R],$$ we can ensure that $\overline{\phi}\leq 0$ on $\R$.
Indeed, assume by contradiction that there exists $x_0\in
\R\setminus[-R,R]$ such that
$$\overline{\phi}(x_0)=\sup_{\R}\overline{\phi}>0.$$Then
$$\left\{%
\begin{array}{ll}
    \I[\overline{\phi},x_0]\leq 0; \\
    \I[\overline{\phi},x_0]-W''(\phi(x_0))\overline{\phi}(x_0)\geq
0; \\
W''(\phi(x_0))>0,\\
\end{array}%
\right.
$$from which $$\overline{\phi}(x_0)\leq 0,$$ a
contradiction. Therefore $\overline{\phi}\leq 0$ on $\R$ and then,
by renaming the constants, from \eqref{phi'infinity} we get
$\phi''\leq \frac{K_1}{1+x^2}$.

 To prove that $\phi''\geq
-\frac{K_1}{1+x^2}$, we look at the infimum of the function
$\phi''+C\phi_a'$ to get similarly that $\phi''+C\phi_a'\geq 0$ on
$\R$.

To show \eqref{phi'''infinity} we proceed as in the proof of
\eqref{phi''infinity}. Indeed, the function $\phi'''$ which is
bounded and of class $C^{1,\beta}$, satisfies
$$\I[\phi''']=W''(\phi)\phi'''+3W'''(\phi)\phi'\phi''+W^{IV}(\phi)(\phi')^3=W''(\phi)\phi'''+o\left(\frac{1}{1+x^2}\right),$$
as $|x|\rightarrow+\infty$, by \eqref{phi'infinity} and
\eqref{phi''infinity}. Then, as before, for $C$ and $a$ large
enough $\phi'''-C\phi_a'\leq 0$ and $\phi'''+C\phi_a'\geq 0$ on
$\R$, which implies \eqref{phi'''infinity}.\finedim

\subsection{Proof of Lemma \ref{psiinfinitylem}.}
Let us prove \eqref{psiinfinity}.

 For $a>0$ we denote by
$\phi_a(x)=\phi\left(\frac{x}{a}\right)$, which is solution of
$$\I[\phi_a]=\frac{1}{a}W'(\phi_a)\quad\text{in }\R.$$ Let $a$ and $b$ be positive numbers, then making a Taylor
expansion of the derivatives of $W$, we get
\begin{equation*}\begin{split}
\I[\psi-(\phi_a-\phi_b)]&=W''(\phi)\psi+\frac{L}{\al}(W''(\phi)-W''(0))+c\phi'+\left(\frac{1}{b}W'(\phi_b)-\frac{1}{a}W'(\phi_a)\right)\\&=
W''(\phi)(\psi-(\phi_a-\phi_b))+W''(\widetilde{\phi})(\phi_a-\phi_b)+\frac{L}{\al}(W''(\widetilde{\phi})-W''(0))\\&+c\phi'
+\left(\frac{1}{b}W'(\widetilde{\phi}_b)-\frac{1}{a}W'(\widetilde{\phi}_a)\right)\\&
=W''(\phi)(\psi-(\phi_a-\phi_b))+W''(0)(\phi_a-\phi_b)+\frac{L}{\al}W'''(0)\widetilde{\phi}
+c\phi'\\&+W''(0)\left(\frac{1}{b}\widetilde{\phi}_b-\frac{1}{a}\widetilde{\phi}_a\right)
+(\phi_a-\phi_b)O(\widetilde{\phi})+O(\widetilde{\phi})^2+O(\widetilde{\phi}_a)^2+
O(\widetilde{\phi}_b)^2,
\end{split}\end{equation*}and then the function
$\overline{\psi}=\psi-(\phi_a-\phi_b)$ satisfies
\begin{equation*}\begin{split}\I[\overline{\psi}]-W''(\phi)\overline{\psi}&=\al(\phi_a-\phi_b)+\frac{L}{\al}W'''(0)\widetilde{\phi}+c\phi'
+\al\left(\frac{1}{b}\widetilde{\phi}_b-\frac{1}{a}\widetilde{\phi}_a\right)\\&+(\phi_a-\phi_b)O(\widetilde{\phi})+O(\widetilde{\phi})^2
+O(\widetilde{\phi}_a)^2+O(\widetilde{\phi}_b)^2.\end{split}\end{equation*}
We want to estimate the right-hand side of the last equality. By
Lemma \ref{phiinfinitylem}, for $|x|\geq\max\{1,|a|,|b|\}$ we have
$$\al(\phi_a-\phi_b)+\frac{L}{\al}W'''(0)\widetilde{\phi}\geq- \frac{1}{\pi
x}\left[(a-b)+\frac{L}{\al^2}W'''(0)\right]-\frac{K_1\al}{
x^2}\left(a^2+b^2+\frac{|L|}{\al^2}|W'''(0)|\right).$$ Choose
$a,b>0$ such that $(a-b)+\frac{L}{\al^2}W'''(0)=0$, then
$$\al(\phi_a-\phi_b)+\frac{L}{\al}W'''(0)\widetilde{\phi}\geq
-\frac{C}{x^2},$$ for $|x|\geq\max\{1,|a|,|b|\}$. Here and in what
follows, as usual $C$ denotes various positive constants. From
Lemma \ref{phiinfinitylem} we also derive that
$$\al\left(\frac{1}{b}\widetilde{\phi}_b-\frac{1}{a}\widetilde{\phi}_a\right)\geq
-\frac{C}{x^2},$$ $$c\phi'\geq -\frac{C}{1+x^2},$$ and
$$(\phi_a-\phi_b)O(\widetilde{\phi})+O(\widetilde{\phi})^2+O(\widetilde{\phi}_a)^2+O(\widetilde{\phi}_b)^2\geq
-\frac{C}{1+x^2},$$ for $|x|\geq\max\{1,|a|,|b|\}$. Then we
conclude that there exists $R>0$ such that for $|x|\geq R$ we have
$$\I[\overline{\psi}]-W''(\phi)\overline{\psi}\geq
-\frac{C}{1+x^2}.$$Now, let us consider the function
$\phi'_d(x)=\phi'\left(\frac{x}{d}\right)$, $d>0$, which is
solution of
$$\I[\phi_d']=\frac{1}{d}W''(\phi_d)\phi_d'\quad\text{in }\R,$$ and denote
$$\overline{\overline{\psi}}=\overline{\psi}-\widetilde{C}\phi_d',$$with $\widetilde{C}>0$. Then, for $|x|\geq R$ we have
$$\I[\overline{\overline{\psi}}]\geq
W''(\phi)\overline{\psi}-\frac{\widetilde{C}}{d}W''(\phi_d)\phi_d'-\frac{C}{1+x^2}=W''(\phi)\overline{\overline{\psi}}
+\widetilde{C}\phi_d'\left(W''(\phi)-\frac{1}{d}W''(\phi_d)\right)-\frac{C}{1+x^2}.$$
Let us choose $d>0$ and $R_2>R$ such that
\begin{equation*}\left\{%
\begin{array}{ll}W''(\phi)-\frac{1}{d}W''(\phi_d)>\frac{1}{2}W''(0)>0&\text{on }\R\setminus[-R_2,R_2];\\
W''(\phi)>0&\text{on }\R\setminus[-R_2,R_2],\\
\end{array}%
\right.\end{equation*} then from \eqref{phi'infinity}, for
$\widetilde{C}$ large enough we get
\begin{equation*}\I[\overline{\overline{\psi}}]-W''(\phi)\overline{\overline{\psi}}\geq
0\quad\text{on }\R\setminus[-R_2,R_2],\end{equation*} and
$$\overline{\overline{\psi}}<0\quad \text{on
}[-R_2,R_2].$$ As in the proof of Lemma \ref{phiinfinitylem}, we
deduce that $\overline{\overline{\psi}}\leq 0$ on $\R$ and then
$$\psi\leq \frac{K_2}{x}+\frac{K_3}{x^2}\quad\text{for
}|x|\geq1,$$ for some $K_2\in\R$ and $K_3>0$.

Looking at the function $\psi
-(\phi_a-\phi_b)+\widetilde{C}\phi_d'$, we conclude similarly that
$$\psi\geq \frac{K_2}{x}-\frac{K_3}{x^2}\quad\text{for }|x|\geq1,$$ and \eqref{psiinfinity} is proved.

Now let us turn to \eqref{psi'infinity}. By deriving the first
equation in \eqref{psi}, we see that the function $\psi'$  which
is bounded and of class $C^{2,\beta}$, is a solution of
$$\I[\psi']=W''(\phi)\psi'+W'''(\phi)\phi'\psi+\frac{L}{\al}W'''(\phi)\phi'+c\phi''\quad\text{in }\R.$$
 Then the function $\overline{\psi}'=\psi'-C\phi_a'$, satisfies
\begin{equation*}\begin{split}\I[\overline{\psi}']-W''(\phi)\overline{\psi}'&=C\phi'_a\left(W''(\phi)-\frac{1}{a}W''(\phi_a)\right)+W'''(\phi)\phi'\psi
+\frac{L}{\al}W'''(\phi)\phi'+c\phi''\\&=C\phi'_a\left(W''(\phi)-\frac{1}{a}W''(\phi_a)\right)+O\left(\frac{1}{1+x^2}\right),
\end{split}\end{equation*}by \eqref{phi'infinity}, \eqref{phi''infinity} and
\eqref{psiinfinity}, and as in the proof of Lemma
\ref{phiinfinitylem}, we deduce that for $C$ and $a$ large enough
$\overline{\psi}'\leq 0$ on $\R$, which implies that $\psi'\leq
\frac{K_3}{1+x^2}$. The inequality $\psi'\geq -\frac{K_3}{1+x^2}$
is obtained similarly by proving that
$\overline{\psi}'+C\phi'_a\geq 0$ on $\R$.

Finally, with the same proof as before, using
\eqref{phi'infinity}-\eqref{psi'infinity}, we can prove the
estimate \eqref{psi''infinity} for the function $\psi''$ which is
a bounded $C^{1,\beta}$ solution of
\begin{equation*}\begin{split}\I[\psi'']&=W''(\phi)\psi''+2W'''(\phi)\phi'\psi'+W^{IV}(\phi)(\phi')^2\psi
+W'''(\phi)\phi''\psi+\frac{L}{\al}W'''(\phi)\phi''
\\&+\frac{L}{\al}W^{IV}(\phi)(\phi')^2+c\phi'''\\&=W''(\phi)\psi''+O\left(\frac{1}{1+x^2}\right).\end{split}\end{equation*}
\finedim

\subsection{Proof of Claims 1-5.}
{\mbox{ }} \bigskip

\noindent {\bf Proof of Claim 1.}\\
We have for $n>|i_0|$
\begin{equation*}\begin{split}\sum_{i=-n\atop i\neq
i_0}^n\frac{1}{x-i}&=\sum_{i=-n}^{i_0-1}\frac{1}{i_0+\gamma-i}+\sum_{i=i_0+1}^n\frac{1}{i_0+\gamma-i}=\sum_{i=1}^{
n+i_0}\frac{1}{i+\gamma}-\sum_{i=1}^{
n-i_0}\frac{1}{i-\gamma}\\&=\left\{%
\begin{array}{ll}
    \sum_{i=1}^{
n}\frac{-2\gamma}{i^2-\gamma^2}, & \hbox{if } i_0=0\\
    \sum_{i=1}^{
n-i_0}\frac{-2\gamma}{i^2-\gamma^2}+\sum_{i=n-i_0+1}^{
n+i_0}\frac{1}{i+\gamma}, & \hbox{if }i_0>0 \\
\sum_{i=1}^{
n+i_0}\frac{-2\gamma}{i^2-\gamma^2}-\sum_{i=n+i_0+1}^{
n-i_0}\frac{1}{i-\gamma}, & \hbox{if }i_0<0 \\
\end{array}%
\right.\rightarrow -2\gamma\sum_{i=1}^{
+\infty}\frac{1}{i^2-\gamma^2}\quad\text{as }n\rightarrow+\infty.
\end{split}\end{equation*}Let us prove the second limit of the
claim.
\begin{equation*}
\sum_{i=-n}^{i_0-1}\frac{1}{(x-i)^2}=
\sum_{i=1}^{n+i_0}\frac{1}{(i+\gamma)^2}\rightarrow
\sum_{i=1}^{+\infty}\frac{1}{(i+\gamma)^2}\quad\text{as
}n\rightarrow+\infty.\end{equation*} Finally
\begin{equation*}\sum_{i=i_0+1}^{n}\frac{1}{(x-i)^2}=
\sum_{i=1}^{n-i_0}\frac{1}{(i-\gamma)^2}\rightarrow
\sum_{i=1}^{+\infty}\frac{1}{(i-\gamma)^2}\quad\text{as
}n\rightarrow+\infty,
\end{equation*}
and the claim  is proved.

By Claim 1 $\sum_{i=-n\atop i\neq i_0}^n\frac{1}{x-i}$,
$\sum_{i=-n}^{i_0-1}\frac{1}{(x-i)^2}$ and
$\sum_{i=i_0+1}^{n}\frac{1}{(x-i)^2}$ are Cauchy sequences and
then for $k>m>|i_0|$ we have
\begin{equation}\label{cauchy1}
\sum_{i=-k}^{-m-1}\frac{1}{x-i}+\sum_{i=m+1}^{k}\frac{1}{x-i}\rightarrow0\quad\text{as
}m,k\rightarrow+\infty,
\end{equation}
\begin{equation}\label{cauchy2}\sum_{i=-k}^{-m-1}\frac{1}{(x-i)^2}
\rightarrow0\quad\text{as }m,k\rightarrow+\infty,\end{equation}
and
\begin{equation}\label{cauchy3}\sum_{i=m+1}^{k}\frac{1}{(x-i)^2}\rightarrow0\quad\text{as
}m,k\rightarrow+\infty.\end{equation}

\noindent {\bf Proof of Claim 2.}\\
We show that $\{s_{\delta,n}^L(x)\}_n$ is a Cauchy sequence. Fix
$x\in\R$ and let $i_0\in\Z$ be the closest integer to $x$ such
that $x=i_0+\gamma$, with
$\gamma\in\left(-\frac{1}{2},\frac{1}{2}\right]$ and $|x-i|\geq
\frac{1}{2}$ for $i\neq i_0$. Let $\delta$ be so small that
$\frac{1}{\delta|p_0|}\geq 2$, then $\frac{|x-i|}{\delta|p_0|}\geq
1$ for $i\neq i_0$. Let $k> m>|i_0|$, using  \eqref{phiinfinity}
and
 \eqref{psiinfinity} we get
\begin{equation*}\begin{split}s_{\delta,k}^L(x)-s_{\delta,m}^L(x)&=
-(k-m)+\sum_{i=-k}^{-m-1}\left[\phi(x_i)+\delta\psi(x_i)\right]
+\sum_{i=m+1}^{k}\left[\phi(x_i)+\delta\psi(x_i)\right] \\&
=\sum_{i=-k}^{-m-1}[\left(\phi(x_i)-1\right)+\delta\psi(x_i)]
+\sum_{i=m+1}^{k}\left[\phi(x_i)+\delta\psi(x_i)\right]\\& \leq
-\left(\frac{1}{\al\pi}-\delta
K_2\right)\delta|p_0|\sum_{i=-k}^{-m-1}\frac{1}{x-i}+(K_1+\delta
K_3)\delta^2|p_0|^2\sum_{i=-k}^{-m-1}\frac{1}{(x-i)^2}\\&-\left(\frac{1}{\al\pi}-\delta
K_2\right)\delta|p_0|\sum_{i=m+1}^{k}\frac{1}{x-i}+(K_1+\delta
K_3)\delta^2|p_0|^2\sum_{i=m+1}^{k}\frac{1}{(x-i)^2},
\end{split}\end{equation*}
and
\begin{equation*}\begin{split}
s_{\delta,k}^L(x)-s_{\delta,m}^L(x)&\geq
-\left(\frac{1}{\al\pi}-\delta
K_2\right)\delta|p_0|\left(\sum_{i=-k}^{-m-1}\frac{1}{x-i}+
\sum_{i=m+1}^{k}\frac{1}{x-i}\right)\\& -(K_1+\delta
K_3)\delta^2|p_0|^2\left(\sum_{i=-k}^{-m-1}\frac{1}{(x-i)^2}+\sum_{i=m+1}^{k}\frac{1}{(x-i)^2}\right).
\end{split}\end{equation*}
Then from \eqref{cauchy1}, \eqref{cauchy2}, \eqref{cauchy3}, we
conclude that
$$|s_{\delta,k}^L(x)-s_{\delta,m}^L(x)|\rightarrow0 \quad\text{as
}m,k\rightarrow+\infty,$$ as desired.\\

\noindent {\bf Proof of Claim 3.}\\
To prove the
uniform convergence, it suffices to show that
$\{(s_{\delta,n}^L)'(x)\}_n$ is a Cauchy sequence uniformly on
compact sets. Let us consider a bounded interval $[a,b]$ and let
$x\in[a,b]$. For $\frac{1}{\delta|p_0|}\geq 2$ and $k>
m>1/2+\max\{|a|,|b|\}$, by \eqref{phi'infinity} and
\eqref{psi'infinity} we have
\begin{equation*}\begin{split}(s_{\delta,k}^L)'(x)-(s_{\delta,m}^L)'(x)&=
\frac{1}{\delta|p_0|}\sum_{i=-k}^{-m-1}\left[\phi'(x_i)+\delta\psi'(x_i)\right]
+\frac{1}{\delta|p_0|}\sum_{i=m+1}^{k}\left[\phi'(x_i)+\delta\psi'(x_i)\right]\\&
\leq(K_1+\delta
K_3)\delta|p_0|\left[\sum_{i=-k}^{-m-1}\frac{1}{(x-i)^2}+\sum_{i=m+1}^{k}\frac{1}{(x-i)^2}\right]\\&
\leq(K_1+\delta
K_3)\delta|p_0|\left[\sum_{i=-k}^{-m-1}\frac{1}{(a-i)^2}+\sum_{i=m+1}^{k}\frac{1}{(b-i)^2}\right],
\end{split}\end{equation*}and
\begin{equation*}(s_{\delta,k}^L)'(x)-(s_{\delta,m}^L)'(x)\geq
-K_3\delta^2|p_0|\left[\sum_{i=-k}^{-m-1}\frac{1}{(a-i)^2}+\sum_{i=m+1}^{k}\frac{1}{(b-i)^2}\right].\end{equation*}
Then by \eqref{cauchy2} and \eqref{cauchy3}
$$\sup_{x\in[a,b]}|(s_{\delta,k}^L)'(x)-(s_{\delta,m}^L)'(x)|\rightarrow0
\quad\text{as }k,m\rightarrow+\infty,$$and Claim 3 is proved.\\

\noindent {\bf Proof of Claim 4.}\\
Claim 4 can be proved like Claim 3. Indeed
$$(s_{\delta,n}^L)''(x)=\frac{1}{\delta^2|p_0|^2}\sum_{i=-n}^{n}\left[\phi''(x_i)
+\delta\psi''(x_i)\right]$$ and using \eqref{phi''infinity} and
\eqref{psi''infinity}, it is easy to show that
$\{(s_{\delta,n}^L)''\}_n$ is a Cauchy sequence uniformly on
compact sets.\\

\noindent {\bf Proof of Claim 5.}\\
We have
$$\I[\phi]=W'(\phi)=W'(\widetilde{\phi})=W''(0)\widetilde{\phi}+O(\widetilde{\phi})^2.$$
Let $x=i_0+\gamma$ with
$\gamma\in\left(-\frac{1}{2},\frac{1}{2}\right]$, and $k>
m>|i_0|$. From \eqref{phiinfinity}, \eqref{cauchy1},
\eqref{cauchy2} and \eqref{cauchy3} we get
\begin{equation*}\begin{split}&\sum_{i=-k}^k\I[\phi,x_i]-\sum_{i=-m}^m\I[\phi,x_i]
=\sum_{i=-k}^{-m-1}[\al\widetilde{\phi}(x_i)+O(\widetilde{\phi}(x_i))^2]
+\sum_{i=m+1}^k[\al\widetilde{\phi}(x_i)+O(\widetilde{\phi}(x_i))^2]\\&
\leq
-\frac{\delta|p_0|}{\pi}\left[\sum_{i=-k}^{-m-1}\frac{1}{x-i}+\sum_{i=m+1}^k\frac{1}{x-i}\right]
+C\sum_{i=-k}^{-m-1}\frac{1}{(x-i)^2}+C\sum_{i=m+1}^k\frac{1}{(x-i)^2}\rightarrow
0,\end{split}\end{equation*}as $m,k\rightarrow+\infty$, for some
constant $C>0$, and
\begin{equation*}\begin{split}&\sum_{i=-k}^k\I[\phi,x_i]-\sum_{i=-m}^m\I[\phi,x_i]
\\&\geq-\frac{\delta|p_0|}{\pi}\left[\sum_{i=-k}^{-m-1}\frac{1}{x-i}+\sum_{i=m+1}^k\frac{1}{x-i}\right]
-C\sum_{i=-k}^{-m-1}\frac{1}{(x-i)^2}-C\sum_{i=m+1}^k\frac{1}{(x-i)^2}\rightarrow
0,\end{split}\end{equation*} as $m,k\rightarrow+\infty$. Then
$\sum_{i=-n}^n\I[\phi,x_i]$ is a Cauchy sequence, i.e. it
converges.

Let us consider now  $\sum_{i=-n}^n\I[\psi,x_i]$.
By \eqref{i1psi}, \eqref{phiinfinity}, \eqref{phi'infinity} and
\eqref{psiinfinity} we get
\begin{equation*}\begin{split}&\sum_{i=-k}^k\I[\psi,x_i]-\sum_{i=-m}^m\I[\psi,x_i]\\&
\leq
\widetilde{C}\left[\sum_{i=-k}^{-m-1}\frac{1}{x-i}+\sum_{i=m+1}^k\frac{1}{x-i}\right]
+C\sum_{i=-k}^{-m-1}\frac{1}{(x-i)^2}+C\sum_{i=m+1}^k\frac{1}{(x-i)^2},
\end{split}\end{equation*}and
\begin{equation*}\begin{split}&\sum_{i=-k}^k\I[\psi,x_i]-\sum_{i=-m}^m\I[\psi,x_i]\\&
\geq
\widetilde{C}\left[\sum_{i=-k}^{-m-1}\frac{1}{x-i}+\sum_{i=m+1}^k\frac{1}{x-i}\right]
-C\sum_{i=-k}^{-m-1}\frac{1}{(x-i)^2}-C\sum_{i=m+1}^k\frac{1}{(x-i)^2},
\end{split}\end{equation*}
 for some $\widetilde{C}\in\R$ and $C>0$, which ensures the convergence of
$\sum_{i=-n}^n\I[\psi,x_i]$.\\

\noindent {\bf Proof of Lemma \ref{lem::s105}}\\
Let $\frac{1}{\delta|p_0|}\geq2$. We first remark that if $z>n+\frac{1}{2}$, then $z_i=\frac{z-i}{\delta|p_0|}\geq 1$ for $i=-n,...,n$ and by \eqref{phiinfinity} and
\eqref{psiinfinity} we have
\begin{equation*}\begin{split}s_{\delta,n}^L(z)&=\frac{L\delta}{\al}+n+1+\sum_{i=-n}^{n}[\phi(z_i)-1
+\delta\psi(z_i)]\\&\leq
\frac{L\delta}{\al}+n+1+\sum_{i=-n}^{n}\left[-\left(\frac{1}{\al\pi}-\delta
K_2\right)\frac{\delta|p_0|}{z-i}+(K_1+\delta
K_3)\frac{\delta^2|p_0|^2}{(z-i)^2}\right],
\end{split}\end{equation*} and
\begin{equation*}\begin{split}s_{\delta,n}^L(z)&\geq\frac{L\delta}{\al}+n+1+\sum_{i=-n}^{n}\left[-\left(\frac{1}{\al\pi}-\delta
K_2\right)\frac{\delta|p_0|}{z-i}-(K_1+\delta
K_3)\frac{\delta^2|p_0|^2}{(z-i)^2}\right].
\end{split}\end{equation*}
By Claim 1, the quantities $\sum_{i=-n}^{n}\frac{1}{z-i}$ and
$\sum_{i=-n}^{n}\frac{1}{(z-i)^2}$ are uniformly bounded on $\R$
by a constant independent of $n$. Hence, we get
\beq\label{s(n)y>n}n-C\leq s_{\delta,n}^L(z)\leq n+C\quad\text{if
}z>n+\frac{1}{2}.\eeq The same argument shows that

\beq\label{s(n)y<ny=n}-n-C\leq s_{\delta,n}^L(z)\leq
-n+C\quad\text{if }z<-n-\frac{1}{2}.\eeq 

If  $|z|<n-\frac{1}{2}$,  then $n>|j_0|$, where $j_0$ is the closest integer to $z$, and  as in the proof of (ii) of Proposition \ref{hproperties} (see Step 1 there), we get 

\begin{equation*}\begin{split}s_{\delta,n}^L(z)-z\leq
\frac{L\delta}{\al}+\frac{3}{2}+\delta\|\psi\|_\infty+\sum_{i=-n\atop
i\neq j_0}^{n}\left[-\left(\frac{1}{\al\pi}-\delta
K_2\right)\frac{\delta|p_0|}{z-i}+(K_1+\delta
K_3)\frac{\delta^2|p_0|^2}{(z-i)^2}\right],
\end{split}\end{equation*} and 
\begin{equation*}\begin{split}s_{\delta,n}^L(z)-z\geq
\frac{L\delta}{\al}-\frac{1}{2}-\delta\|\psi\|_\infty+\sum_{i=-n\atop
i\neq j_0}^{n}\left[-\left(\frac{1}{\al\pi}-\delta
K_2\right)\frac{\delta|p_0|}{z-i}-(K_1+\delta
K_3)\frac{\delta^2|p_0|^2}{(z-i)^2}\right].
\end{split}\end{equation*} 
Then, again by Claim 1

 \beq\label{sn(y)-y} -C\leq s_{\delta,n}^L(z)-z\leq C
\quad\text{if }|z|<n-\frac{1}{2}.\eeq

Now, let $i_0\in\Z$ be the closest integer to $x$, let us assume
$n>|i_0|+1+a.$ We have
\begin{equation*}\begin{split}\int_{|y|\geq
a}[s_{\delta,n}^L(x+y)-s_{\delta,n}^L(x)]\mu(dy)&=\int_{a\leq|y|<
n-1-|i_0|}[s_{\delta,n}^L(x+y)-s_{\delta,n}^L(x)]\mu(dy)\\&+\int_{n-1-|i_0|\leq|y|\leq
n+1+|i_0|}[...]\mu(dy)+\int_{|y|>
n+1+|i_0|}[...]\mu(dy).\end{split}\end{equation*}

If $|y|< n-1-|i_0|$, then $|x+y|<n-\frac{1}{2}$ and by \eqref{sn(y)-y}
\begin{equation*}\begin{split}\int_{a\leq|y|<
n-1-|i_0|}[s_{\delta,n}^L(x+y)-s_{\delta,n}^L(x)]\mu(dy)&\leq\int_{a\leq|y|\leq
n-1-|i_0|}(y+2C)\mu(dy)\\&=\int_{a\leq|y|\leq
n-1-|i_0|}2C\mu(dy)\leq \frac{2C}{a},
\end{split}\end{equation*}and
$$\int_{a\leq|y|<
n-1-|i_0|}[s_{\delta,n}^L(x+y)-s_{\delta,n}^L(x)]\mu(dy)\ge
-\frac{2C}{a}.$$ Then
\beq\label{limalimn1}\lim_{a\rightarrow+\infty}\lim_{n\rightarrow+\infty}\int_{a\leq|y|\leq
n-1-|i_0|}[s_{\delta,n}^L(x+y)-s_{\delta,n}^L(x)]\mu(dy)=0.\eeq

Next, since $|s_{\delta,n}^L(z)|\leq Cn$ for any $z\in\R$, we have
\begin{equation}\label{seconinteni0}\begin{split}&\left|\int_{n-1-|i_0|\leq|y|\leq
n+1+|i_0|}[s_{\delta,n}^L(x+y)-s_{\delta,n}^L(x)]\mu(dy)\right|\\&\leq
Cn\int_{n-1-|i_0|\leq|y|\leq
n+1+|i_0|}\mu(dy)=\widetilde{C}\frac{n(|i_0|+1)}{n^2-(|i_0|+1)^2}\rightarrow0
\quad\text{as }n\rightarrow+\infty.\end{split}\end{equation}

Finally, if $y>n+1+|i_0|$, then $x+y>n+\frac{1}{2}$, while if  $y<-n-1-|i_0|$, then
$x+y<-n-\frac{1}{2}$. Hence, using \eqref{s(n)y>n} and \eqref{s(n)y<ny=n}, we
obtain
\begin{equation*}\begin{split}&\int_{|y|>
n+1+|i_0|}[s_{\delta,n}^L(x+y)-s_{\delta,n}^L(x)]\mu(dy)\\&=\int_{y>
n+1+|i_0|}[s_{\delta,n}^L(x+y)-s_{\delta,n}^L(x)]\mu(dy)
+\int_{y<-n-1-|i_0|}[s_{\delta,n}^L(x+y)-s_{\delta,n}^L(x)]\mu(dy)\\&
\leq \int_{y> n+1+|i_0|}[n+C-s_{\delta,n}^L(x)]\mu(dy)+\int_{y<
-n-1-|i_0|}[-n+C-s_{\delta,n}^L(x)]\mu(dy)\\& =\int_{|y|>
n+1+|i_0|}[C-s_{\delta,n}^L(x)]\mu(dy),
\end{split}\end{equation*}
and
$$\int_{|y|>
n+1+|i_0|}[s_{\delta,n}^L(x+y)-s_{\delta,n}^L(x)]\mu(dy)\geq
\int_{|y|> n+1+|i_0|}[-C-s_{\delta,n}^L(x)]\mu(dy).$$ We deduce
that
$$\lim_{n\rightarrow+\infty}\int_{|y|>
n+1+|i_0|}[s_{\delta,n}^L(x+y)-s_{\delta,n}^L(x)]\mu(dy)=0.$$
Hence, by the previous limit, \eqref{limalimn1} and
\eqref{seconinteni0}, we derive \eqref{limalimnI2}.
This ends the proof of Lemma \ref{lem::s105}.


\begin{thebibliography}{10}













\bibitem {bkm}
{\scshape P. Biler, G. Karch, R. Monneau}  
Nonlinear diffusion of dislocation density and self-similar solutions,
\emph{Communications in Mathematical Physics}, {\bf 294} (2010), no. 1,
145-168. 




\bibitem{csm}{\sc X. Cabr\'{e} and J. Sol\`{a}-Morales}, Layer
solutions in a half-space for boundary reactions, {\em Comm. Pure
Appl. Math.}, {\bf 58} (2005) no. 12, 1678-1732.








\bibitem{dmopt}{\sc A.~De~Masi, E.~Orlandi, E.~Presutti  and L.~Triolo},  Motion by
  curvature by scaling nonlocal evolution equations, {\em J. Statist. Phys.}, {\bf 73}
  (1993),  543--570.

\bibitem{Denoual}
{\sc C. Denoual},
Dynamic dislocation modeling by combining Peierls Nabarro and Galerkin methods,
{\em Phys. Rev. B}, {\bf 70} (2004), 024106.

\bibitem{di}{\sc J. Droniou and C. Imbert }, Fractal first order partial differential
equations, {\em Archive for Rational Mechanics and Analysis}, {\bf
182} (2006), no. 2, 299-331.










\bibitem{fim2}{\sc N. Forcadel, C. Imbert and R. Monneau}, Homogenization of fully overdamped Frenkel-Kontorova models,
{\em Journal of Differential Equations}, {\bf 246} (2009), no. 1,
1057-1097.











\bibitem{gonzalezmonneau}{\sc M. Gonz\'{a}lez and R. Monneau},
Slow motion of particle systems as a limit of a reaction-diffusion
equation with half-Laplacian in dimension one, preprint
hal-00497492.


\bibitem{H} {\sc A. K. Head}
Dislocation group dynamics III. 
Similarity solutions of the continuum approximation,
{\em Phil. Magazine}, {\bf 26}, (1972), 65-72.


\bibitem{hl}{\sc J. R. Hirth and L. Lothe}, Theory of
dislocations, Second Edition. Malabar, Florida: Krieger, 1992.


















\bibitem{lpv} {\sc P. L. Lions, G. C. Papanicolaou and S. R. S. Varadhan}, Homogenization of Hamilton-Jacobi equations, unpublished, 1986.


\bibitem{mp1}{\sc R. Monneau and S. Patrizi},  Homogenization of the Peierls-Nabarro model for dislocation dynamics, {\em preprint}.

\bibitem{MBW}
{\sc A.B. Movchan, R. Bullough, J.R. Willis}, Stability of a
dislocation: discrete model, {\em Eur. J. Appl. Math.} {\bf 9}
(1998), 373-396.







\bibitem{N}
{\sc F.R.N. Nabarro},
Fifty-year study of the Peierls-Nabarro stress,
{\em Material Science and Engineering A} {\bf 234-236} (1997), 67-76.











\bibitem{WXM}
{\sc H. Wei, Y. Xiang, P. Ming},
A Generalized Peierls-Nabarro Model for Curved Dislocations
Using Discrete Fourier Transform,
{\em Communications in computational physics} {\bf 4}(2) (2008), 275-293.




\end{thebibliography}
\end{document}